\documentclass[conference, 10pt]{IEEEtran}

\usepackage{cite}

\usepackage[cmex10]{amsmath}
\usepackage{url}
\usepackage{graphicx}
\usepackage{color}
\usepackage{placeins}
\usepackage{float}
\usepackage{tabularx,colortbl}
\usepackage{ifthen}
\usepackage{amssymb,amsfonts}
\usepackage{algorithmic}
\usepackage{textcomp}
\usepackage{comment}
\usepackage{xcolor}
\usepackage{siunitx}
\usepackage{bm}
\usepackage{tikz}
\usetikzlibrary{shapes.geometric}
\usepackage{pgfplots}
\pgfplotsset{compat=1.5}
\usepgfplotslibrary{units}
\pdfmapfile{-mpfonts.map}
\newcommand{\NsTD}
{(0.045836623610465865, 7.312514667796357e17)
 (0.004583662361046586, 8.714932660767088e13)
 (0.00045836623610465875, 8.696162763524142e11)
 (4.583662361046588e-5, 8.696212078859776e9)
 (4.583662361046588e-6, 8.696209659241791e7)
 (4.5836623610465875e-7, 869622.0414815104)
 (4.583662361046588e-8, 8698.050159155375)
 (4.583662361046588e-9, 88.85228050676926)}
 
 \newcommand{\NsReg}
{ (0.045836623610465865, 89.888646819922)
 (0.004583662361046586, 89.88864681992617)
 (0.00045836623610465875, 89.88864682030312)
 (4.583662361046588e-5, 89.88864685805079)
 (4.583662361046588e-6, 89.88865063274206)
 (4.5836623610465875e-7, 89.88902810284226)
 (4.583662361046588e-8, 89.92678513447186)
 (4.583662361046588e-9, 93.80477851053323)}
 
  \newcommand{\NsCal}
{(0.045836623610465865, 2.2280187502799396)
 (0.004583662361046586, 2.2280187502802216)
 (0.00045836623610465875, 2.228018750308169)
 (4.583662361046588e-5, 2.2280187531026403)
 (4.583662361046588e-6, 2.2280190325498253)
 (4.5836623610465875e-7, 2.228046976255404)
 (4.583662361046588e-8, 2.2308351237564854)
 (4.583662361046588e-9, 2.4878710299433755)}

 \newcommand{\dtreg}
{(1/0.2911, 43.81961506660352)
 (1/0.19809884900797972, 93.80477851053323)
 (1/0.14964520636475284, 159.55523309671983)
 (1/0.12011532126570548, 246.38289256069515)}
  \newcommand{\dtTD}
{(1/0.2911, 40.1489964920729)
 (1/0.19809884900797972, 88.85228050676926)
 (1/0.14964520636475284, 152.3268701262547)
 (1/0.12011532126570548, 236.33097016862362)}
  \newcommand{\dtcal}
{ (1/0.2911, 2.392697592780229)
 (1/0.19809884900797972, 2.4878710299433755)
 (1/0.14964520636475284, 2.52613869675355)
 (1/0.12011532126570548, 2.5448140390352854)}
 
 \newcommand{\pointCal}{ (-3.246760839074665e-5, 9.616929373485522e-18)
 (-2.38732414637843e-5, 1.0262179003323002e-10)
 (-1.5278874536821953e-5, 6.866899007011487e-6)
 (-6.684507609859608e-6, 0.002877184565318088)
 (1.9098593171027432e-6, 0.007546164445041387)
 (1.0504226244065088e-5, 0.0001241217032641909)
 (1.9098593171027432e-5, 1.2865959251189453e-8)
 (2.7692960097989777e-5, 8.468104206940436e-15)
 (3.6287327024952135e-5, 3.5674020549787556e-23)
 (4.488169395191448e-5, 9.556778996679139e-34)
 (5.3476060878876824e-5, 2.356112348715686e-45)
 (6.207042780583917e-5, 2.688079856569271e-51)
 (7.066479473280151e-5, 2.693313755199618e-57)}
 
 \newcommand{\pointReg}{(-3.533239736640076e-5, 1.3848383543085732e-20)
 (-2.6738030439438415e-5, 8.101807302153243e-13)
 (-1.8143663512476067e-5, 2.9615527923697107e-7)
 (-9.549296585513723e-6, 0.0006754609843162455)
 (-9.549296585513716e-7, 0.00961505840494168)
 (7.639437268410973e-6, 0.0008574415635203621)
 (1.6233804195373317e-5, 4.8351149572229e-7)
 (2.4828171122335662e-5, 1.752036746372305e-12)
 (3.342253804929802e-5, 4.165574295795849e-20)
 (4.2016904976260365e-5, 6.62469066968374e-30)
 (5.061127190322271e-5, 7.038082512112376e-42)
 (5.9205638830185054e-5, 2.61089543617288e-49)
 (6.78000057571474e-5, 1.859138790149327e-55)}
 
 \newcommand{\pointTD}{ (-3.819718634205488e-5, 7.091538600199686e-24)
 (-2.960281941509253e-5, 3.6800620597588426e-15)
 (-2.1008452488130182e-5, 7.2922616994538655e-9)
 (-1.2414085561167838e-5, 9.023820572942651e-5)
 (-3.8197186342054865e-6, 0.006972583934476463)
 (4.774648292756858e-6, 0.0033708416163701757)
 (1.3369015219719203e-5, 1.0250611577509501e-5)
 (2.1963382146681554e-5, 1.1059198268452464e-10)
 (3.0557749073643905e-5, 1.2137944814022047e-10)
 (3.915211600060625e-5, 1.5539755030152144e-10)
 (4.7746482927568594e-5, 1.8941562773764797e-10)
 (5.634084985453094e-5, 2.2343370500403626e-10)
 (6.493521678149328e-5, 2.57451782072662e-10)
 (7.352958370845563e-5, 2.9146985891543674e-10)}

\newcommand{\CurrentCal}{ 
 (-3.819718634205488e-5, 1.2172375252489374e-23)
 (-3.762422854692405e-5, 5.111896165966853e-23)
 (-3.7051270751793234e-5, 2.1592547278956307e-22)
 (-3.647831295666241e-5, 8.904808682566061e-22)
 (-3.590535516153159e-5, 3.590709977848416e-21)
 (-3.533239736640076e-5, 1.4156418818411788e-20)
 (-3.475943957126994e-5, 5.456835701258916e-20)
 (-3.418648177613912e-5, 2.056568643431504e-19)
 (-3.361352398100829e-5, 7.578095289547406e-19)
 (-3.304056618587747e-5, 2.7301795742243967e-18)
 (-3.246760839074665e-5, 9.616929373485522e-18)
 (-3.189465059561582e-5, 3.312039179734809e-17)
 (-3.1321692800485e-5, 1.1152393697325913e-16)
 (-3.074873500535418e-5, 3.671589065686471e-16)
 (-3.0175777210223353e-5, 1.1818249450828371e-15)
 (-2.960281941509253e-5, 3.7193333657066425e-15)
 (-2.9029861619961708e-5, 1.1444311969592424e-14)
 (-2.8456903824830886e-5, 3.4429178276621486e-14)
 (-2.7883946029700063e-5, 1.0126882561061036e-13)
 (-2.7310988234569238e-5, 2.9123045978636097e-13)
  (-2.6738030439438415e-5, 8.188596742049274e-13)
 (-2.6165072644307593e-5, 2.2510938397014795e-12)
 (-2.5592114849176767e-5, 6.050467094345105e-12)
 (-2.5019157054045945e-5, 1.5899928044633666e-11)
 (-2.4446199258915123e-5, 4.08518878392795e-11)
 (-2.38732414637843e-5, 1.0262179003323002e-10)
 (-2.3300283668653478e-5, 2.52044539255915e-10)
 (-2.2727325873522653e-5, 6.052363765247934e-10)
 (-2.215436807839183e-5, 1.4209622818239531e-9)
 (-2.1581410283261008e-5, 3.2617423533638175e-9)
 (-2.1008452488130182e-5, 7.320253354019735e-9)
 (-2.043549469299936e-5, 1.6062447214908096e-8)
 (-1.9862536897868538e-5, 3.4459283768198736e-8)
 (-1.9289579102737715e-5, 7.227856238742973e-8)
 (-1.8716621307606893e-5, 1.4822494297160946e-7)
 (-1.8143663512476067e-5, 2.971949396344643e-7)
 (-1.7570705717345245e-5, 5.825989175613947e-7)
 (-1.6997747922214423e-5, 1.1166212791980552e-6)
 (-1.6424790127083597e-5, 2.0924250749691096e-6)
 (-1.5851832331952775e-5, 3.833553735737771e-6)
 (-1.5278874536821953e-5, 6.866899007011487e-6)
 (-1.470591674169113e-5, 1.2026162882237608e-5)
 (-1.4132958946560308e-5, 2.0592101959522475e-5)
 (-1.3560001151429482e-5, 3.447317827949199e-5)
 (-1.298704335629866e-5, 5.642465322752332e-5)
 (-1.2414085561167838e-5, 9.029491465862435e-5)
 (-1.1841127766037012e-5, 0.00014127470480645519)
 (-1.126816997090619e-5, 0.0002161086570775647)
 (-1.0695212175775368e-5, 0.0003232112650381238)
 (-1.0122254380644545e-5, 0.00047261480747037496)
 (-9.549296585513723e-6, 0.0006756699857066075)
 (-8.976338790382897e-6, 0.0009444269457164907)
 (-8.403380995252075e-6, 0.001290650156064799)
 (-7.830423200121253e-6, 0.0017244679879438235)
 (-7.257465404990427e-6, 0.0022527250345286266)
 (-6.684507609859608e-6, 0.002877184565318088)
 (-6.111549814728782e-6, 0.0035928062882045354)
 (-5.538592019597957e-6, 0.004386382509516166)
 (-4.965634224467138e-6, 0.005235834533006952)
 (-4.392676429336312e-6, 0.006110435339861306)
  (-3.8197186342054865e-6, 0.006972128059821423)
 (-3.2467608390746676e-6, 0.007777959505916156)
 (-2.673803043943842e-6, 0.008483465517985418)
 (-2.100845248813023e-6, 0.009046663166392289)
 (-1.5278874536821973e-6, 0.009432162355574661)
 (-9.549296585513716e-7, 0.009614840609335817)
 (-3.819718634205527e-7, 0.009582551417640736)
 (1.9098593171027297e-7, 0.009337460411022692)
 (7.639437268410919e-7, 0.008895805088466175)
 (1.3369015219719176e-6, 0.008286115141012712)
 (1.9098593171027432e-6, 0.007546164445041387)
 (2.482817112233562e-6, 0.006719107131122905)
 (3.055774907364388e-6, 0.0058493457912746705)
 (3.6287327024952135e-6, 0.004978676569012334)
 (4.201690497626032e-6, 0.004143161912041022)
 (4.774648292756858e-6, 0.003371023534789867)
 (5.347606087887684e-6, 0.0026816623408258177)
 (5.920563883018503e-6, 0.002085736311033276)
 (6.493521678149328e-6, 0.0015860917639753185)
 (7.066479473280147e-6, 0.0011792650972366956)
  (7.639437268410973e-6, 0.000857253886388407)
 (8.212395063541799e-6, 0.0006092885934080326)
 (8.785352858672618e-6, 0.0004234020071477855)
 (9.358310653803443e-6, 0.00028767346034408984)
 (9.931268448934262e-6, 0.0001911016647654043)
 (1.0504226244065088e-5, 0.0001241217032641909)
 (1.1077184039195913e-5, 7.88225331176041e-5)
 (1.1650141834326732e-5, 4.8941052120961916e-5)
 (1.2223099629457558e-5, 2.9711004912351796e-5)
 (1.2796057424588384e-5, 1.7635322439528704e-5)
 (1.3369015219719203e-5, 1.0234636285296725e-5)
 (1.3941973014850028e-5, 5.807447766448263e-6)
 (1.4514930809980854e-5, 3.22198168682674e-6)
 (1.5087888605111673e-5, 1.7477796023677118e-6)
 (1.56608464002425e-5, 9.269944113827028e-7)
 (1.6233804195373317e-5, 4.807235970905917e-7)
 (1.6806761990504143e-5, 2.4374891034348246e-7)
 (1.737971978563497e-5, 1.208425552292217e-7)
 (1.7952677580765788e-5, 5.857713990555873e-8)
 (1.8525635375896613e-5, 2.7763138237149525e-8)
  (1.9098593171027432e-5, 1.2865959251189453e-8)
 (1.9671550966158258e-5, 5.829752452345595e-9)
 (2.0244508761289084e-5, 2.582816139685259e-9)
 (2.0817466556419903e-5, 1.118854001056662e-9)
 (2.1390424351550728e-5, 4.739046977745895e-10)
 (2.1963382146681554e-5, 1.9626704919531077e-10)
 (2.2536339941812373e-5, 7.947740399682128e-11)
 (2.31092977369432e-5, 3.146885682695772e-11)
 (2.3682255532074024e-5, 1.2183172116618827e-11)
 (2.425521332720485e-5, 4.6119332260335795e-12)
 (2.4828171122335662e-5, 1.707064562903155e-12)
 (2.5401128917466488e-5, 6.178204255195944e-13)
 (2.5974086712597313e-5, 2.186356687180588e-13)
 (2.654704450772814e-5, 7.565319009695252e-14)
 (2.7120002302858965e-5, 2.559658182664139e-14)
 (2.7692960097989777e-5, 8.468104206940436e-15)
 (2.8265917893120603e-5, 2.739308066054112e-15)
 (2.8838875688251428e-5, 8.664568939539441e-16)
 (2.9411833483382254e-5, 2.679818356191734e-16)
 (2.998479127851308e-5, 8.104326652399634e-17)
  (3.0557749073643905e-5, 2.3965281484223385e-17)
 (3.113070686877472e-5, 6.929516592415785e-18)
 (3.170366466390554e-5, 1.959200546588559e-18)
 (3.227662245903637e-5, 5.416404431059254e-19)
 (3.2849580254167194e-5, 1.4641991095310733e-19)
 (3.342253804929802e-5, 3.8703176155078244e-20)
 (3.399549584442883e-5, 1.0003476463429822e-20)
 (3.456845363955966e-5, 2.5282123077633224e-21)
 (3.5141411434690484e-5, 6.247907496877752e-22)
 (3.571436922982131e-5, 1.50978169539852e-22)
 (3.6287327024952135e-5, 3.5674020549787556e-23)
 (3.686028482008295e-5, 8.242288510994206e-24)
 (3.743324261521377e-5, 1.8620893136355315e-24)
 (3.80062004103446e-5, 4.1134729089229592e-25)
 (3.8579158205475424e-5, 8.885254605248662e-26)
 (3.915211600060625e-5, 1.8766387968907178e-26)
 (3.972507379573706e-5, 3.875572498248365e-27)
 (4.029803159086789e-5, 7.825781937736648e-28)
 (4.087098938599871e-5, 1.5450593143849836e-28)
 (4.144394718112954e-5, 2.982440945693061e-29)
  (4.2016904976260365e-5, 5.6283904369374045e-30)
 (4.258986277139119e-5, 1.0383600760161155e-30)
 (4.3162820566522e-5, 1.8724626601796244e-31)
 (4.373577836165283e-5, 3.299951717614589e-32)
 (4.4308736156783654e-5, 5.682306236010486e-33)
 (4.488169395191448e-5, 9.556778996679139e-34)
 (4.5454651747045305e-5, 1.5690556474833888e-34)
 (4.602760954217612e-5, 2.512892718788359e-35)
 (4.660056733730694e-5, 3.921235538768368e-36)
 (4.717352513243777e-5, 5.951109747464909e-37)
 (4.7746482927568594e-5, 8.764237456545663e-38)
 (4.831944072269942e-5, 1.2456272813901008e-38)
 (4.8892398517830246e-5, 1.7036276275714145e-39)
 (4.946535631296106e-5, 2.1921947622634116e-40)
 (5.0038314108091884e-5, 2.658626007861849e-41)
 (5.061127190322271e-5, 2.7765429253711403e-42)
 (5.1184229698353535e-5, 1.92399354958654e-43)
 (5.175718749348436e-5, 2.1978181159257892e-44)
 (5.233014528861517e-5, 1.6581413070365605e-44)
 (5.2903103083746e-5, 3.717503937940918e-45)
  (5.3476060878876824e-5, 2.356112348715686e-45)
 (5.404901867400765e-5, 5.908070235670459e-46)
 (5.4621976469138475e-5, 1.747021089122086e-46)
 (5.519493426426929e-5, 1.2173249310431744e-47)
 (5.576789205940011e-5, 1.5617362088640235e-47)
 (5.634084985453094e-5, 1.4374331204512865e-47)
 (5.6913807649661765e-5, 8.123310058718623e-48)
 (5.748676544479259e-5, 3.7572332786245574e-48)
 (5.80597232399234e-5, 1.485786806371678e-48)
 (5.863268103505423e-5, 4.995579675276375e-49)
 (5.9205638830185054e-5, 1.3129735858496017e-49)
 (5.977859662531588e-5, 1.6646485133500163e-50)
 (6.0351554420446705e-5, 9.036091588049133e-51)
 (6.092451221557753e-5, 9.512723276020202e-51)
 (6.149747001070834e-5, 5.659058097478209e-51)
 (6.207042780583917e-5, 2.688079856569271e-51)
 (6.264338560097e-5, 1.0898305061723892e-51)
 (6.321634339610082e-5, 3.77285598859902e-52)
 (6.378930119123165e-5, 1.044965293438366e-52)
 (6.436225898636246e-5, 1.657186320421867e-53)
  (6.493521678149328e-5, 4.696134933425465e-54)
 (6.550817457662411e-5, 6.242764387689377e-54)
 (6.608113237175493e-5, 3.916100537601695e-54)
 (6.665409016688576e-5, 1.915938016609599e-54)
 (6.722704796201659e-5, 7.96034156905874e-55)
 (6.78000057571474e-5, 2.833647016838409e-55)
 (6.837296355227822e-5, 8.217233927804374e-56)
 (6.894592134740905e-5, 1.5215717269164424e-56)
 (6.951887914253988e-5, 2.053268257529176e-57)
 (7.00918369376707e-5, 4.036471817211925e-57)
 (7.066479473280151e-5, 2.693313755199618e-57)
 (7.123775252793234e-5, 1.3598402402806467e-57)
 (7.181071032306316e-5, 5.791131026249704e-58)
 (7.238366811819399e-5, 2.1172281519401293e-58)
 (7.295662591332482e-5, 6.396562133935033e-59)
 (7.352958370845563e-5, 1.32969306714403e-59)
 (7.410254150358645e-5, 5.079706921340269e-61)
 (7.467549929871728e-5, 2.5618196283007635e-60)
 (7.52484570938481e-5, 1.8397559135066322e-60)
 (7.582141488897893e-5, 9.609296775925087e-61)
}

\newcommand{\CurrentTD}{
(-3.819718634205488e-5, 7.091538600199686e-24)
 (-3.762422854692405e-5, 4.2383547419553845e-23)
 (-3.7051270751793234e-5, 2.0109144108928815e-22)
 (-3.647831295666241e-5, 8.601916893398586e-22)
 (-3.590535516153159e-5, 3.509460936202973e-21)
 (-3.533239736640076e-5, 1.389243617222868e-20)
 (-3.475943957126994e-5, 5.3644261350936636e-20)
 (-3.418648177613912e-5, 2.0238431672604314e-19)
 (-3.361352398100829e-5, 7.463597046087027e-19)
 (-3.304056618587747e-5, 2.6909098570806983e-18)
 (-3.246760839074665e-5, 9.485232173928551e-18)
 (-3.189465059561582e-5, 3.26888710497392e-17)
 (-3.1321692800485e-5, 1.1014288972354917e-16)
 (-3.074873500535418e-5, 3.628422900534351e-16)
 (-3.0175777210223353e-5, 1.168649222698929e-15)
 (-2.960281941509253e-5, 3.6800620597588426e-15)
 (-2.9029861619961708e-5, 1.1330018622374352e-14)
 (-2.8456903824830886e-5, 3.410440138191253e-14)
 (-2.7883946029700063e-5, 1.0036779375811515e-13)
 (-2.7310988234569238e-5, 2.8879008606120745e-13)
 (-2.6738030439438415e-5, 8.124075632870099e-13)
 (-2.6165072644307593e-5, 2.234442766771271e-12)
 (-2.5592114849176767e-5, 6.0085259387470545e-12)
 (-2.5019157054045945e-5, 1.5796828653518913e-11)
 (-2.4446199258915123e-5, 4.060457557235811e-11)
 (-2.38732414637843e-5, 1.020429460859074e-10)
 (-2.3300283668653478e-5, 2.507227820949663e-10)
 (-2.2727325873522653e-5, 6.02292230862606e-10)
 (-2.215436807839183e-5, 1.4145660846695906e-9)
 (-2.1581410283261008e-5, 3.2481913836433274e-9)
 (-2.1008452488130182e-5, 7.2922616994538655e-9)
 (-2.043549469299936e-5, 1.6006081945722884e-8)
 (-1.9862536897868538e-5, 3.4348666054031337e-8)
 (-1.9289579102737715e-5, 7.206703747026277e-8)
 (-1.8716621307606893e-5, 1.4783093756282973e-7)
 (-1.8143663512476067e-5, 2.9648026155637916e-7)
 (-1.7570705717345245e-5, 5.813370006140457e-7)
 (-1.6997747922214423e-5, 1.1144531556898693e-6)
 (-1.6424790127083597e-5, 2.0888021250968766e-6)
 (-1.5851832331952775e-5, 3.8276690649185244e-6)
 (-1.5278874536821953e-5, 6.857614149143287e-6)
 (-1.470591674169113e-5, 1.201194365114708e-5)
 (-1.4132958946560308e-5, 2.0570986284925578e-5)
 (-1.3560001151429482e-5, 3.444280837530285e-5)
 (-1.298704335629866e-5, 5.638241176006785e-5)
 (-1.2414085561167838e-5, 9.023820572942651e-5)
 (-1.1841127766037012e-5, 0.00014120141079789432)
 (-1.126816997090619e-5, 0.00021601777744967329)
 (-1.0695212175775368e-5, 0.0003231037010623277)
 (-1.0122254380644545e-5, 0.0004724942002107685)
 (-9.549296585513723e-6, 0.0006755434521636127)
 (-8.976338790382897e-6, 0.0009443055054750822)
 (-8.403380995252075e-6, 0.0012905486183824428)
 (-7.830423200121253e-6, 0.0017244040974513989)
 (-7.257465404990427e-6, 0.0022527177663242397)
 (-6.684507609859608e-6, 0.0028772516166410567)
 (-6.111549814728782e-6, 0.0035929610137033136)
 (-5.538592019597957e-6, 0.0043866307497076994)
 (-4.965634224467138e-6, 0.005236171936083932)
 (-4.392676429336312e-6, 0.00611084582032963)
 (-3.8197186342054865e-6, 0.006972583934476463)
 (-3.2467608390746676e-6, 0.00777842359330576)
 (-2.673803043943842e-6, 0.008483895117741301)
 (-2.100845248813023e-6, 0.009047015426007398)
 (-1.5278874536821973e-6, 0.009432400182966564)
 (-9.549296585513716e-7, 0.009614938095619295)
 (-3.819718634205527e-7, 0.009582497759111394)
 (1.9098593171027297e-7, 0.009337261573016559)
 (7.639437268410919e-7, 0.008895482883548207)
 (1.3369015219719176e-6, 0.008285703903703692)
  (1.9098593171027432e-6, 0.007545705949856963)
 (2.482817112233562e-6, 0.0067186447134359945)
 (3.055774907364388e-6, 0.005848918723358052)
 (3.6287327024952135e-6, 0.0049783156087574035)
 (4.201690497626032e-6, 0.004142886629029434)
 (4.774648292756858e-6, 0.0033708416163701757)
 (5.347606087887684e-6, 0.002681570690703563)
 (5.920563883018503e-6, 0.0020857234593326234)
 (6.493521678149328e-6, 0.001586140959559516)
 (7.066479473280147e-6, 0.0011793574566426235)
 (7.639437268410973e-6, 0.000857371098198373)
 (8.212395063541799e-6, 0.0006094148701873344)
 (8.785352858672618e-6, 0.0004235251825204892)
 (9.358310653803443e-6, 0.0002877853007637637)
 (9.931268448934262e-6, 0.00019119757753941456)
 (1.0504226244065088e-5, 0.00012420006845034339)
 (1.1077184039195913e-5, 7.888387738727055e-5)
 (1.1650141834326732e-5, 4.898723673437335e-5)
 (1.2223099629457558e-5, 2.9744538649463627e-5)
 (1.2796057424588384e-5, 1.7658850819539035e-5)
 (1.3369015219719203e-5, 1.0250611577509501e-5)
 (1.3941973014850028e-5, 5.81795438764334e-6)
 (1.4514930809980854e-5, 3.228677750990974e-6)
 (1.5087888605111673e-5, 1.7519138299682459e-6)
 (1.56608464002425e-5, 9.294635617026541e-7)
 (1.6233804195373317e-5, 4.821447783599018e-7)
 (1.6806761990504143e-5, 2.44530465798402e-7)
 (1.737971978563497e-5, 1.212449696260968e-7)
 (1.7952677580765788e-5, 5.876091710383906e-8)
 (1.8525635375896613e-5, 2.7823954867845993e-8)
 (1.9098593171027432e-5, 1.2859059026226276e-8)
 (1.9671550966158258e-5, 5.786063808956456e-9)
 (2.0244508761289084e-5, 2.5191425266311162e-9)
 (2.0817466556419903e-5, 1.0440611304759187e-9)
 (2.1390424351550728e-5, 3.9253791998457273e-10)
 (2.1963382146681554e-5, 1.1059198268452464e-10)
 (2.2536339941812373e-5, 9.407529231957672e-12)
 (2.31092977369432e-5, 6.010775497753198e-11)
 (2.3682255532074024e-5, 8.184716380984086e-11)
 (2.425521332720485e-5, 9.176573128368854e-11)
  (2.4828171122335662e-5, 9.697159486434097e-11)
 (2.5401128917466488e-5, 1.0034217672848539e-10)
 (2.5974086712597313e-5, 1.0301457287587948e-10)
 (2.654704450772814e-5, 1.0542749255331599e-10)
 (2.7120002302858965e-5, 1.0774620039746136e-10)
 (2.7692960097989777e-5, 1.1003148768690491e-10)
 (2.8265917893120603e-5, 1.1230519142897892e-10)
 (2.8838875688251428e-5, 1.145749722726573e-10)
 (2.9411833483382254e-5, 1.168434549604303e-10)
 (2.998479127851308e-5, 1.1911151787299945e-10)
 (3.0557749073643905e-5, 1.2137944814022047e-10)
 (3.113070686877472e-5, 1.236473374466762e-10)
 (3.170366466390554e-5, 1.2591521439183769e-10)
 (3.227662245903637e-5, 1.2818308769102583e-10)
 (3.2849580254167194e-5, 1.304509599389573e-10)
 (3.342253804929802e-5, 1.3271883189036518e-10)
 (3.399549584442883e-5, 1.3498670375973885e-10)
 (3.456845363955966e-5, 1.3725457560663087e-10)
 (3.5141411434690484e-5, 1.3952244744719165e-10)
 (3.571436922982131e-5, 1.4179031928570143e-10)
 (3.6287327024952135e-5, 1.440581911232613e-10)
 (3.686028482008295e-5, 1.4632606296014272e-10)
 (3.743324261521377e-5, 1.4859393479640622e-10)
 (3.80062004103446e-5, 1.508618066320606e-10)
 (3.8579158205475424e-5, 1.5312967846710217e-10)
 (3.915211600060625e-5, 1.5539755030152144e-10)
 (3.972507379573706e-5, 1.576654221353123e-10)
 (4.029803159086789e-5, 1.5993329396846724e-10)
 (4.087098938599871e-5, 1.622011658009747e-10)
 (4.144394718112954e-5, 1.6446903763282745e-10)
 (4.2016904976260365e-5, 1.6673690946401751e-10)
 (4.258986277139119e-5, 1.690047812945368e-10)
 (4.3162820566522e-5, 1.7127265312437756e-10)
 (4.373577836165283e-5, 1.7354052495353162e-10)
 (4.4308736156783654e-5, 1.7580839678198992e-10)
 (4.488169395191448e-5, 1.7807626860974428e-10)
 (4.5454651747045305e-5, 1.8034414043678495e-10)
 (4.602760954217612e-5, 1.8261201226310382e-10)
 (4.660056733730694e-5, 1.848798840886928e-10)
 (4.717352513243777e-5, 1.871477559135433e-10)
  (4.7746482927568594e-5, 1.8941562773764797e-10)
 (4.831944072269942e-5, 1.9168349956099998e-10)
 (4.8892398517830246e-5, 1.9395137138358814e-10)
 (4.946535631296106e-5, 1.9621924320540503e-10)
 (5.0038314108091884e-5, 1.9848711502644137e-10)
 (5.061127190322271e-5, 2.0075498684668956e-10)
 (5.1184229698353535e-5, 2.030228586661431e-10)
 (5.175718749348436e-5, 2.0529073048479126e-10)
 (5.233014528861517e-5, 2.075586023026253e-10)
 (5.2903103083746e-5, 2.0982647411964052e-10)
 (5.3476060878876824e-5, 2.12094345935827e-10)
 (5.404901867400765e-5, 2.143622177511761e-10)
 (5.4621976469138475e-5, 2.166300895656796e-10)
 (5.519493426426929e-5, 2.1889796137932938e-10)
 (5.576789205940011e-5, 2.2116583319211649e-10)
 (5.634084985453094e-5, 2.2343370500403626e-10)
 (5.6913807649661765e-5, 2.2570157681507957e-10)
 (5.748676544479259e-5, 2.2796944862523264e-10)
 (5.80597232399234e-5, 2.3023732043449112e-10)
 (5.863268103505423e-5, 2.325051922428439e-10)
  (5.9205638830185054e-5, 2.3477306405028464e-10)
 (5.977859662531588e-5, 2.37040935856802e-10)
 (6.0351554420446705e-5, 2.3930880766238836e-10)
 (6.092451221557753e-5, 2.4157667946703983e-10)
 (6.149747001070834e-5, 2.4384455127074454e-10)
 (6.207042780583917e-5, 2.4611242307349515e-10)
 (6.264338560097e-5, 2.483802948752843e-10)
 (6.321634339610082e-5, 2.506481666761035e-10)
 (6.378930119123165e-5, 2.529160384759463e-10)
 (6.436225898636246e-5, 2.5518391027480024e-10)
 (6.493521678149328e-5, 2.57451782072662e-10)
 (6.550817457662411e-5, 2.597196538695214e-10)
 (6.608113237175493e-5, 2.619875256653688e-10)
 (6.665409016688576e-5, 2.642553974601969e-10)
 (6.722704796201659e-5, 2.665232692539952e-10)
 (6.78000057571474e-5, 2.6879114104675645e-10)
 (6.837296355227822e-5, 2.7105901283847195e-10)
 (6.894592134740905e-5, 2.73326884629133e-10)
 (6.951887914253988e-5, 2.7559475641873067e-10)
 (7.00918369376707e-5, 2.7786262820725734e-10)
  (7.066479473280151e-5, 2.801304999947073e-10)
 (7.123775252793234e-5, 2.823983717810645e-10)
 (7.181071032306316e-5, 2.846662435663236e-10)
 (7.238366811819399e-5, 2.8693411535047884e-10)
 (7.295662591332482e-5, 2.892019871335178e-10)
 (7.352958370845563e-5, 2.9146985891543674e-10)
 (7.410254150358645e-5, 2.937377306962264e-10)
 (7.467549929871728e-5, 2.96005602475877e-10)
 (7.52484570938481e-5, 2.982734742543818e-10)
 (7.582141488897893e-5, 3.005413460317301e-10)
 }
 
 \newcommand{\CurrentReg}{
 (-3.819718634205488e-5, 1.157288008887032e-23)
 (-3.762422854692405e-5, 4.971719139163085e-23)
 (-3.7051270751793234e-5, 2.107318110069309e-22)
 (-3.647831295666241e-5, 8.694714208991156e-22)
 (-3.590535516153159e-5, 3.509434383729139e-21)
 (-3.533239736640076e-5, 1.3848383543085732e-20)
 (-3.475943957126994e-5, 5.342824569049517e-20)
 (-3.418648177613912e-5, 2.0153453891461872e-19)
 (-3.361352398100829e-5, 7.432497324939056e-19)
 (-3.304056618587747e-5, 2.6799485376394925e-18)
 (-3.246760839074665e-5, 9.447660831742201e-18)
 (-3.189465059561582e-5, 3.256326992927407e-17)
 (-3.1321692800485e-5, 1.0973301983520798e-16)
 (-3.074873500535418e-5, 3.615363680366796e-16)
 (-3.0175777210223353e-5, 1.1645864198181825e-15)
 (-2.960281941509253e-5, 3.667720713318409e-15)
 (-2.9029861619961708e-5, 1.1293416287383838e-14)
 (-2.8456903824830886e-5, 3.399841712810981e-14)
 (-2.7883946029700063e-5, 1.0006819886491716e-13)
 (-2.7310988234569238e-5, 2.8796336322378843e-13)
  (-2.6738030439438415e-5, 8.101807302153243e-13)
 (-2.6165072644307593e-5, 2.2285883325660197e-12)
 (-2.5592114849176767e-5, 5.993504329815644e-12)
 (-2.5019157054045945e-5, 1.575921512857098e-11)
 (-2.4446199258915123e-5, 4.051267330098401e-11)
 (-2.38732414637843e-5, 1.0182385837596054e-10)
 (-2.3300283668653478e-5, 2.5021325157509334e-10)
 (-2.2727325873522653e-5, 6.011363074599923e-10)
 (-2.215436807839183e-5, 1.4120084670722068e-9)
 (-2.1581410283261008e-5, 3.2426728339235343e-9)
 (-2.1008452488130182e-5, 7.280651881902652e-9)
 (-2.043549469299936e-5, 1.5982272055812438e-8)
 (-1.9862536897868538e-5, 3.4301074229533343e-8)
 (-1.9289579102737715e-5, 7.197434356751179e-8)
 (-1.8716621307606893e-5, 1.4765506328297473e-7)
 (-1.8143663512476067e-5, 2.9615527923697107e-7)
 (-1.7570705717345245e-5, 5.807523723086609e-7)
 (-1.6997747922214423e-5, 1.1134296314008895e-6)
 (-1.6424790127083597e-5, 2.0870590154218325e-6)
 (-1.5851832331952775e-5, 3.824782786419787e-6)
  (-1.5278874536821953e-5, 6.852970296618309e-6)
 (-1.470591674169113e-5, 1.200468868040803e-5)
 (-1.4132958946560308e-5, 2.0559990227560496e-5)
 (-1.3560001151429482e-5, 3.4426656371180076e-5)
 (-1.298704335629866e-5, 5.6359448332158786e-5)
 (-1.2414085561167838e-5, 9.020665947595953e-5)
 (-1.1841127766037012e-5, 0.0001411596255800623)
 (-1.126816997090619e-5, 0.00021596456598885775)
 (-1.0695212175775368e-5, 0.00032303881565131124)
 (-1.0122254380644545e-5, 0.0004724188826493716)
 (-9.549296585513723e-6, 0.0006754609843162455)
 (-8.976338790382897e-6, 0.000944221640783278)
 (-8.403380995252075e-6, 0.001290471745307142)
 (-7.830423200121253e-6, 0.0017243449824741253)
 (-7.257465404990427e-6, 0.002252688759529646)
 (-6.684507609859608e-6, 0.002877265283112298)
 (-6.111549814728782e-6, 0.003593028273659423)
 (-5.538592019597957e-6, 0.004386758749988845)
 (-4.965634224467138e-6, 0.0052363620030921425)
 (-4.392676429336312e-6, 0.00611109192125484)
 (-3.8197186342054865e-6, 0.00697287206507911)
 (-3.2467608390746676e-6, 0.007778732398357092)
 (-2.673803043943842e-6, 0.008484197841896246)
 (-2.100845248813023e-6, 0.00904728302491257)
 (-1.5278874536821973e-6, 0.009432605156358275)
 (-9.549296585513716e-7, 0.00961505840494168)
 (-3.819718634205527e-7, 0.009582520123128253)
 (1.9098593171027297e-7, 0.009337183489861244)
 (7.639437268410919e-7, 0.008895313045644277)
 (1.3369015219719176e-6, 0.00828546083074482)
 (1.9098593171027432e-6, 0.007545415193779877)
 (2.482817112233562e-6, 0.006718335127636956)
 (3.055774907364388e-6, 0.005848618507184288)
 (3.6287327024952135e-6, 0.0049780488111303214)
 (4.201690497626032e-6, 0.004142670649246614)
 (4.774648292756858e-6, 0.003370685994834425)
 (5.347606087887684e-6, 0.002681477205099632)
 (5.920563883018503e-6, 0.0020856873018258192)
 (6.493521678149328e-6, 0.001586152613246672)
 (7.066479473280147e-6, 0.0011794048363672819)
  (7.639437268410973e-6, 0.0008574415635203621)
 (8.212395063541799e-6, 0.0006094968281738654)
 (8.785352858672618e-6, 0.00042360913676685145)
 (9.358310653803443e-6, 0.00028786433469360213)
 (9.931268448934262e-6, 0.00019126735902374109)
 (1.0504226244065088e-5, 0.00012425851771827702)
 (1.1077184039195913e-5, 7.893065263166071e-5)
 (1.1650141834326732e-5, 4.9023172176259505e-5)
 (1.2223099629457558e-5, 2.977113228815124e-5)
 (1.2796057424588384e-5, 1.7677857040539235e-5)
 (1.3369015219719203e-5, 1.026375699933585e-5)
 (1.3941973014850028e-5, 5.826768981489486e-6)
 (1.4514930809980854e-5, 3.234418451776088e-6)
 (1.5087888605111673e-5, 1.7555528577148094e-6)
 (1.56608464002425e-5, 9.317154495998232e-7)
 (1.6233804195373317e-5, 4.8351149572229e-7)
 (1.6806761990504143e-5, 2.453504845271021e-7)
 (1.737971978563497e-5, 1.217380053255484e-7)
 (1.7952677580765788e-5, 5.906458703949135e-8)
 (1.8525635375896613e-5, 2.8021573400725597e-8)
  (1.9098593171027432e-5, 1.2999447184123978e-8)
 (1.9671550966158258e-5, 5.896943622854832e-9)
 (2.0244508761289084e-5, 2.6157818997484114e-9)
 (2.0817466556419903e-5, 1.134622299981217e-9)
 (2.1390424351550728e-5, 4.812592874671987e-10)
 (2.1963382146681554e-5, 1.9961251785875647e-10)
 (2.2536339941812373e-5, 8.096178905559688e-11)
 (2.31092977369432e-5, 3.2111378000067415e-11)
 (2.3682255532074024e-5, 1.2454525022716145e-11)
 (2.425521332720485e-5, 4.723757976000746e-12)
 (2.4828171122335662e-5, 1.752036746372305e-12)
 (2.5401128917466488e-5, 6.354724733969372e-13)
 (2.5974086712597313e-5, 2.2539853122242732e-13)
 (2.654704450772814e-5, 7.818241064374678e-14)
 (2.7120002302858965e-5, 2.651999754076459e-14)
 (2.7692960097989777e-5, 8.797254027521528e-15)
 (2.8265917893120603e-5, 2.8538606938211445e-15)
 (2.8838875688251428e-5, 9.053846094964674e-16)
 (2.9411833483382254e-5, 2.808994059940599e-16)
 (2.998479127851308e-5, 8.522922459073901e-17)
  (3.0557749073643905e-5, 2.5290000215911637e-17)
 (3.113070686877472e-5, 7.338955037479613e-18)
 (3.170366466390554e-5, 2.0827981142756803e-18)
 (3.227662245903637e-5, 5.780828285905126e-19)
 (3.2849580254167194e-5, 1.5691527605764158e-19)
 (3.342253804929802e-5, 4.165574295795849e-20)
 (3.399549584442883e-5, 1.081486761760015e-20)
 (3.456845363955966e-5, 2.7460362561380033e-21)
 (3.5141411434690484e-5, 6.819178264449165e-22)
 (3.571436922982131e-5, 1.6561528695599762e-22)
 (3.6287327024952135e-5, 3.933808175640537e-23)
 (3.686028482008295e-5, 9.138438224988724e-24)
 (3.743324261521377e-5, 2.0762438918278747e-24)
 (3.80062004103446e-5, 4.613534171779127e-25)
 (3.8579158205475424e-5, 1.0026268788981783e-25)
 (3.915211600060625e-5, 2.1310626716011894e-26)
 (3.972507379573706e-5, 4.430014966637447e-27)
 (4.029803159086789e-5, 9.006722413087654e-28)
 (4.087098938599871e-5, 1.790940244747179e-28)
 (4.144394718112954e-5, 3.4829541844449013e-29)
  (4.2016904976260365e-5, 6.62469066968374e-30)
 (4.258986277139119e-5, 1.2323456275697018e-30)
 (4.3162820566522e-5, 2.2420473270614033e-31)
 (4.373577836165283e-5, 3.989299985376987e-32)
 (4.4308736156783654e-5, 6.941917648844464e-33)
 (4.488169395191448e-5, 1.1813504087218933e-33)
 (4.5454651747045305e-5, 1.9659754903632913e-34)
 (4.602760954217612e-5, 3.199162152379259e-35)
 (4.660056733730694e-5, 5.090245495808781e-36)
 (4.717352513243777e-5, 7.916805606425341e-37)
 (4.7746482927568594e-5, 1.2036697132006e-37)
 (4.831944072269942e-5, 1.7876413823977357e-38)
 (4.8892398517830246e-5, 2.5917420990899617e-39)
 (4.946535631296106e-5, 3.679922525238566e-40)
 (5.0038314108091884e-5, 5.009462957810546e-41)
 (5.061127190322271e-5, 7.038082512112376e-42)
 (5.1184229698353535e-5, 7.919725455584004e-43)
 (5.175718749348436e-5, 1.4181983020274282e-43)
 (5.233014528861517e-5, 5.916057637351775e-45)
 (5.2903103083746e-5, 2.203322597024554e-45)
  (5.3476060878876824e-5, 8.968230409945525e-46)
 (5.404901867400765e-5, 6.457741805550531e-46)
 (5.4621976469138475e-5, 3.7639974345408864e-46)
 (5.519493426426929e-5, 1.7349377810564296e-46)
 (5.576789205940011e-5, 6.862909040504049e-47)
 (5.634084985453094e-5, 2.3093449368415023e-47)
 (5.6913807649661765e-5, 6.071854580914246e-48)
 (5.748676544479259e-5, 7.732144115969125e-49)
 (5.80597232399234e-5, 4.155378923460436e-49)
 (5.863268103505423e-5, 4.386042875427158e-49)
 (5.9205638830185054e-5, 2.61089543617288e-49)
 (5.977859662531588e-5, 1.240562461428025e-49)
 (6.0351554420446705e-5, 5.030765582354661e-50)
 (6.092451221557753e-5, 1.7420324839289712e-50)
 (6.149747001070834e-5, 4.827181192237439e-51)
 (6.207042780583917e-5, 7.670387767347766e-52)
 (6.264338560097e-5, 2.157176306386316e-52)
 (6.321634339610082e-5, 2.8767457422347967e-52)
 (6.378930119123165e-5, 1.805556421286607e-52)
 (6.436225898636246e-5, 8.835617690369315e-53)
  (6.493521678149328e-5, 3.671554637381405e-53)
 (6.550817457662411e-5, 1.3071576634165902e-53)
 (6.608113237175493e-5, 3.79158328253757e-54)
 (6.665409016688576e-5, 7.027391226646404e-55)
 (6.722704796201659e-5, 9.418767534590797e-56)
 (6.78000057571474e-5, 1.859138790149327e-55)
 (6.837296355227822e-5, 1.2410128184039289e-55)
 (6.894592134740905e-5, 6.266671986971304e-56)
 (6.951887914253988e-5, 2.668948287950434e-56)
 (7.00918369376707e-5, 9.758114021413019e-57)
 (7.066479473280151e-5, 2.9483841024337684e-57)
 (7.123775252793234e-5, 6.131091301508638e-58)
 (7.181071032306316e-5, 2.3215070834865127e-59)
 (7.238366811819399e-5, 1.1794765226733372e-58)
 (7.295662591332482e-5, 8.472332114900952e-59)
 (7.352958370845563e-5, 4.425372813324334e-59)
 (7.410254150358645e-5, 1.932523960828643e-59)
 (7.467549929871728e-5, 7.249551126629966e-60)
 (7.52484570938481e-5, 2.2729085225684333e-60)
 (7.582141488897893e-5, 5.174017291685474e-61)
 }

\usepackage{cancel}
\usepackage[normalem]{ulem}
\usepackage{soul}
\usepackage{fancyhdr}

\fancypagestyle{firstpage}
{
 \footskip 30pt
\fancyfoot[C]{DOI: 10.1109/AP-S/USNC-URSI47032.2022.9887361 | © 2022 IEEE. Personal use of this material is permitted. Permission from IEEE must be obtained for all other uses, in any current or future media, including reprinting/republishing this material for advertising or promotional purposes, creating new collective works, for resale or redistribution to servers or lists, or reuse of any copyrighted component of this work in other works.}
}

\newcommand{\minus}{\scalebox{0.75}[0.9]{$-$}}
\newcommand{\plus}{\centering \scalebox{0.9}[0.9]{$+$}}

\makeatletter

\newcounter{author}
\renewcommand{\author}[2][]{
   \stepcounter{author}
   \@namedef{author@\theauthor}{#2}
   \@namedef{authorlabel@\theauthor}{#1}
}

\newcounter{address}
\newcommand{\address}[2][]{
   \stepcounter{address}
   \@namedef{address@\theaddress}{#2}
   \@namedef{addresslabel@\theaddress}{#1}
}

\newcommand{\alsep}{and}

\def\newmaketitle{\par%
  \begingroup%
  \normalfont%
  \def\thefootnote{}%  the \thanks{} mark type is empty
  \def\footnotemark{}% and kill space from \thanks within author
  \let\@makefnmark\relax% V1.7, must *really* kill footnotemark to remove all \textsuperscript spacing as well.
  \footnotesize%       equal spacing between thanks lines
  \footnotesep 0.7\baselineskip%see global setting of \footnotesep for more info
  \normalsize%
  \twocolumn[\thenewmaketitle\@IEEEaftertitletext]%
  \if@IEEEusingpubid
     \enlargethispage{-\@IEEEpubidpullup}%
  \fi
  \endgroup
  \setcounter{footnote}{0}\let\maketitle\relax\let\@maketitle\relax
  \gdef\@thanks{}%
  \let\thanks\relax}

\def\thenewmaketitle{
  \newpage
  \begin{center}%
    \vskip0.2em{\Huge\@IEEEcompsoconly{\sffamily}\@IEEEcompsocconfonly{\normalfont\normalsize\vskip 2\@IEEEnormalsizeunitybaselineskip
   \bfseries\large}\@title\par}\vskip1.0em\par%
    \vspace{1ex}
    \newcounter{c@author}
    \newcounter{c@tmp}
    \ifthenelse{\value{author}=2}{%
      \newcommand{\liand}{ and }}{%
      \newcommand{\liand}{, and }}
    \ifthenelse{\value{address}<2}{%
      \@nameuse{author@1}%
      \stepcounter{c@author}%
      \whiledo{\value{c@author}<\value{author}}{%
        \setcounter{c@tmp}{\value{author}}%
        \addtocounter{c@tmp}{-\value{c@author}}%
        \ifthenelse{\value{c@tmp}=1}{%
          \renewcommand{\alsep}{\liand}}{\renewcommand{\alsep}{, }}%
        \stepcounter{c@author}\alsep \@nameuse{author@\thec@author}}\\%
    }
    {%Add address references after the author's name
      \@nameuse{author@1}${}^{(\ref{\@nameuse{authorlabel@1}})}$%
      \stepcounter{c@author}%
      \whiledo{\value{c@author}<\value{author}}{%
      \setcounter{c@tmp}{\value{author}}%
      \addtocounter{c@tmp}{-\value{c@author}}%
      \ifthenelse{\value{c@tmp}=1}{%
        \renewcommand{\alsep}{\liand}}{\renewcommand{\alsep}{, }}%
      \stepcounter{c@author}\alsep \@nameuse{author@\thec@author}%
        ${}^{(\ref{\@nameuse{authorlabel@\thec@author}})}$%
      }
    }
    \vspace{0.2ex}

    \ifthenelse{\value{address}>0}{%
      \ifthenelse{\value{address}=1}{
        {\@nameuse{address@1}}
      }
      {%Output the addresses as an enumerated list
        \newcounter{c@address}

        \begin{center}
        \whiledo{\value{c@address}<\value{address}}
        {
          \refstepcounter{c@address}
            ${}^{(\thec@address)}$\,%
              \label{\@nameuse{addresslabel@\thec@address}}%
              \@nameuse{address@\thec@address}\\ %
        }
        \end{center}
      } % end of the address creation ifthenelse block
    }
    {
      \relax
    }
  \end{center}
}

\makeatother

\title{Calder\'on Preconditioners for the TD-EFIE discretized with Convolution Quadratures}

\author[org1]{Pierrick Cordel}
\author[org1]{Alexandre Dély}
\author[org2]{Adrien Merlini}
\author[org1]{Francesco P. Andriulli}

\address[org1]{Politecnico di Torino, 10129 Turin, Italy, https://www.det.polito.it/it/}
\address[org2]{IMT Atlantique, 29238 Brest, France}

\begin{document}
\newmaketitle
%\vspace{5cm}
%\maketitle
%\vspace{1cm}
\thispagestyle{firstpage}
\begin{abstract} 
This work focuses on the preconditioning and DC stabilization of the time domain electric field integral equation discretized in time with the convolution quadrature method. The standard formulation of the equation suffers from severe ill-conditioning for large time steps and refined meshes, in addition to DC instabilities plaguing standard solutions for late time steps. This work addresses all these issues by preconditioning the TD-EFIE operator matrices with a Calder\'on approach. Numerical results will corroborate the theory, showing the practical relevance of the proposed advancements.
\end{abstract}

\section{Introduction}

The time-domain electric field integral equation (TD-EFIE) is a powerful formulation for modelling the electromagnetic radiation and scattering from perfectly electrically conducting (PEC) objects in the time domain. Among the different discretization strategies of this equation, convolution quadrature methods (CQM) are particularly effective and easily derived from frequency domain implementations \cite{lubich1988convolution,wang2011implicit}. Their effectiveness notwithstanding, these formulations suffer from severe ill-conditioning for large time steps and refined meshes. Moreover, the TD-EFIE operator has a static null space which leads, with numerical and machine precision errors, to the emergence of spurious static currents (DC-instabilities), a phenomenon limiting the simulations' late-time precision. In this work, we address all the above-mentioned limitations by designing a suitable time domain Calder\'on preconditioner for the TD-EFIE formulation discretized by the CQM. Differently from the frequency domain where the preconditioning is generally done on the matrix system, it is found that a preconditioning applied before the time discretization can concurrently solve conditioning issues and DC instabilities. Theoretical considerations and numerical studies confirm the effectiveness of the approach together with its practical 
relevance.

\section{Background and Notation}

Consider a PEC object of boundary $\Gamma$ and outpointing normal $\mathbf{\hat{n}}$ excited by an electromagnetic field $(\mathbf{e}^{\text{inc}}, \mathbf{h}^{\text{inc}}) (\mathbf{r},t)$. The incident field induces a current $\mathbf{j}$ on $\Gamma$ which can be computed by solving the TD-EFIE
\begin{equation}
\eta_{0}\mathcal{\bm{T}}(\mathbf{j})(\mathbf{r},t)=\minus \mathbf{\hat{n}}(r) \times \mathbf{e}^{\text{inc}}(\mathbf{r},t)\,,\quad \forall (\mathbf{r},t) \in \Gamma\times \mathbb{R}\,,
 \label{EFIE}
\end{equation}
where $\eta_{0}$ is the permeability of the background. The TD-EFIE operator $\mathcal{\bm{T}}$ includes the contributions of the vector and scalar potentials, respectively denoted $\mathcal{\bm{T}}_s$ and $\mathcal{\bm{T}}_h$ \cite{dely2019large}
\begin{equation}
 \mathcal{\bm{T}}(f)(\mathbf{r},t)=\minus\frac{1}{c_0}\frac{\partial}{\partial t}\mathcal{\bm{T}}_s(f)(\mathbf{r},t) \plus c_0\int_{\minus\infty}^{t} \hspace{-7px} \mathcal{\bm{T}}_h(f)(\mathbf{r},t')dt',
\end{equation}
where $c_0$ is the speed of light in the medium.

In this study, Rao-Wilton-Glisson (RWG) basis functions $(f^{\text{rwg}}_n)_{N_s}$ and their rotated counterparts $(\mathbf{\hat{n}} \times f^{\text{rwg}}_n)_{N_s}$ have been used as source and tests functions for the spatial discretization, where $N_s$ is the number of edges of the mesh. The time discretization is a convolution quadrature using an implicit Runge-Kutta method (here, 2 stages Radau IIA) with a time step $\Delta t$. The resulting discrete marching-on-in-time (MOT) scheme is
\begin{equation}
 \forall i \in \mathbb{N},\quad \bm{\mathrm{T}}_{0}\mathbf{J}_i=\mathbf{E}_i\minus\sum_{j=1}^{i}\bm{\mathrm{T}}_{j}\mathbf{J}_{i\minus j}\,,
 \label{MOT1}
\end{equation}
where $\mathbf{J}$ and  $\mathbf{E}$ are respectively the array of coefficients of the RWG expansion of the current $\mathbf{j}$ and the array of $\minus\eta_0^{\minus1}\mathbf{\hat{n}} \times\mathbf{e}^{\text{inc}}$ tested with rotated RWG, at different time steps, and 
\begin{equation}
    \begin{split}
        [\widetilde{\bm{\mathrm{T}}}(s)]_{m,n}&=\int_\Gamma\mathbf{\hat{n}}\times f^{\text{rwg}}_m\mathcal{L}\left(\mathcal{\bm{T}}\left(f^{\text{rwg}}_n\delta\right)\right)(s)d\Gamma,\\
        \bm{\mathrm{T}}_i&=\mathcal{Z}^{\minus1}\left(z\mapsto \widetilde{\bm{\mathrm{T}}} \left(\mathbf{s}\left(z\right)\right)\right)_i,
    \end{split}
\end{equation}
where $\mathcal{L}$ is the Laplace transform, $\delta$ is the time Dirac delta, $\mathcal{Z}^{\minus1}$ is the inverse $\mathcal{Z}$-transform, and $\mathbf{s}(z)$
is fully determined by the Runge-Kutta method and $\Delta t$ \cite{dely2019large}.

\section{On a Calder\'on preconditioner for the CQM}

Calder\'on preconditioners are based on the Calder\'on identity $\mathcal{\bm{T}}^2= \minus\mathcal{\bm{I}}/4\plus \mathcal{\bm{K}}^2$, where $\mathcal{\bm{K}}$ is a compact operator and $\mathcal{\bm{I}}$ is the identity. The operator $\mathcal{\bm{T}}^2$ is therefore well-conditioned for large time steps and refined meshes. This yields the following preconditioned TD-EFIE
\begin{equation}
 \eta_0\mathcal{\bm{T}}^2 (\mathbf{j})(\mathbf{r},t)= \minus\mathcal{\bm{T}}(\mathbf{\hat{n}} \times \mathbf{e}^{\text{inc}})(\mathbf{r},t).
 \label{Calderon}
\end{equation}

One could think of directly deriving a preconditioner from the above formula \eqref{MOT1}, similarly to what is done in the frequency domain, and only precondition $\bm{\mathrm{T}}_{0}$ at each step of the MOT.
%One could think of directly deriving a preconditioner from the above formula, similarly to what is done in the frequency domain, and preconditioning each time step of the MOT $\bm{T}_{0}$ only in \eqref{MOT1} \todo[inline]{This sentence does not make sense, do you mean only precond $T_0$ in (3)?}. 
Doing this would indeed solve the conditioning problems but the solution currents would remain unaltered and subject to DC instabilities. This has motivated the development of the new approach presented in this work: instead of preconditioning $\bm{\mathrm{T}}_{0}(s)$ only, 
%we Calder\'on preconditioned the entire Laplace domain matrix $\widetilde{\bm{\mathrm{T}}}$ which resulted in a DC-stable scheme at the price of extra matrix multiplications at the right-hand-side. 
we apply a Calder\'on-type preconditioning to the entire time domain, which results in a DC-stable scheme at the price of extra matrix multiplications at the right-hand-side.
%, which is equivalent to applied the Calder\'on preconditioner to $\widetilde{\bm{\mathrm{T}}}$
%we applied the Calder\'on preconditioner
%we Calder\'on conditioned the entire time domain, which is equivalent to apply the matrix preconditioner to $\widetilde{\bm{\mathrm{T}}}$ and results in a DC-stable scheme at the price of extra matrix multiplications at the right-hand-side.
In particular, after discretizing the TD-EFIE operator with the RWG basis functions and the preconditioning operator with the Buffa-Christiansen (BC) functions $(f^{\text{bc}}_n)_{N_s}$, the preconditioning is done with matrices associated to all time steps. By defining
$[\widetilde{\bm{\mathbb{T}}}(s)]_{m,n}=\int_\Gamma\mathbf{\hat{n}}\times f^{\text{bc}}_m\mathcal{L}(\mathcal{\bm{T}}(f^{\text{bc}}_n\delta))(s)d\Gamma,$ $\bm{\mathbb{T}}_i=\mathcal{Z}^{\minus1}(z\mapsto \widetilde{\bm{\mathbb{T}}}(\mathbf{s}(z)))_i$ and after some manipulations, the following MOT Calder\'on preconditioned scheme
\begin{equation}
 [\bm{\mathbb{T}}\bm{G_m}^{\minus1}* \bm{\mathrm{T}}]_{0} \mathbf{J}_{i} = [\bm{\mathbb{T}}\bm{G_m}^{\minus1} *\mathbf{E}]_i
 \minus\sum_{j=1}^{i}[\bm{\mathbb{T}}\bm{G_m}^{\minus1}*\bm{\mathrm{T}}]_{j}\mathbf{J}_{i\minus j}\,,
  \label{e2}
\end{equation}
is obtained, with $\bm{G_m}$ the gram matrix between the BC and rotated RWG functions and $*$ is the convolution product. However, the MOT in \eqref{e2}, as the one in \eqref{MOT1}, involves unbounded number of large terms in the convolutions, leading to a quadratic complexity with the time step, because of the time-integral in the scalar potential contribution of the operators.
To remove this time integral, the preconditioned EFIE operator and the right hand side of formulation \eqref{Calderon} are evaluated by separating the vector and scalar potential contributions
\begin{align}
 \mathcal{\bm{T}}^2&=\frac{1}{c_0^2}\frac{\partial}{\partial t}\mathcal{\bm{T}}_s\frac{\partial}{\partial t}\mathcal{\bm{T}}_s \minus\mathcal{\bm{T}}_h\mathcal{\bm{T}}_s \minus\mathcal{\bm{T}}_s\mathcal{\bm{T}}_h,
 \label{trick}
 \\
 \mathcal{\bm{T}}(\mathbf{\hat{n}} \times \mathbf{e}^{\text{inc}})&=\minus\frac{1}{c_0}\frac{\partial}{\partial t}\mathcal{\bm{T}}_s(\mathbf{\hat{n}} \times \mathbf{e}^{\text{inc}})\plus c_0\mathcal{\bm{T}}_h(\mathbf{\hat{n}} \times \mathbf{e}^{\text{prim}}),
 \label{trick2}
\end{align}
where $\mathbf{e}^{\text{prim}}=\int_{\minus \infty}^{t}\hspace{-4px}\mathbf{e}^{\text{inc}}$, because $\mathcal{\bm{T}}_{h}^2=0$. The CQM discrete versions of the operators $\mathcal{\bm{T}}^2$, $\minus\frac{1}{c_0}\frac{\partial}{\partial t}\mathcal{\bm{T}}_s$, and $c_0\mathcal{\bm{T}}_h$, respectively denoted by the matrix sequences $(\bm{\mathrm{T}}^\text{cal}_i)$, $(\bm{\mathrm{T}}^\alpha_i)$ and $(\bm{\mathrm{T}}^\beta_i)$ converge to zeros.
The sums can therefore be truncated and we denote by $N_\text{conv}$ the last considered term.
By extending the previous notation on $\mathcal{\bm{T}}_s$ and $\mathcal{\bm{T}}_h$, one can check that 
\begin{equation}
 \begin{split}
 & \bm{\mathrm{T}}^\text{cal}\hspace{-2px}=c_0^{\minus2}\hspace{-1px} \mathcal{Z}^{\minus1}\hspace{-2px}\big(\hspace{-1px}\bm{s}^2\widetilde{\bm{\mathbb{T}}}_s \bm{G_m}^{\minus1} \widetilde{\bm{\mathrm{T}}}_s \minus\widetilde{\bm{\mathbb{T}}}_s \bm{G_m}^{\minus1} \widetilde{\bm{\mathrm{T}}}_h  \minus \widetilde{\bm{\mathbb{T}}}_h \bm{G_m}^{\minus1}  \widetilde{\bm{\mathrm{T}}}_s\hspace{-1px}\big)\,,\\
 &\bm{\mathrm{T}}^\alpha\hspace{-2px} = \minus c_0^{\minus1} \hspace{-1px} \mathcal{Z}^{\minus1}\hspace{-2px}\big(\bm{s}\widetilde{\bm{\mathbb{T}}}_s\big) \bm{G_m}^{\minus1} \text{ and }
 \bm{\mathrm{T}}^\beta \hspace{-2px}= c_0 \mathcal{Z}^{\minus1} \hspace{-2px}\big(\widetilde{\bm{\mathbb{T}}}_h\big) \bm{G_m}^{\minus1}\,.
 \end{split}
\end{equation}
The MOT \eqref{e2} is therefore rewritten as
\begin{equation}
 \bm{\mathrm{T}}^\text{cal}_0\mathbf{J}_{i}=\sum_{j=0}^{N_{\text{conv}}}\left(\bm{\mathrm{T}}^\alpha_j\mathbf{E}_{i\minus j}\plus\bm{\mathrm{T}}^\beta_j\mathbf{E}_{i\minus j}^{\text{prim}}\right)
 \minus\sum_{j=1}^{N_{\text{conv}}}\bm{\mathrm{T}}^\text{cal}_{j}\mathbf{J}_{i\minus j}\,,
\end{equation}
where $\mathbf{E}_{i}^{\text{prim}}$ is the array of $\minus\eta_{0}^{\minus1}\mathbf{\hat{n}} \times \mathbf{e}^{\text{prim}}$ tested with rotated RWG at different time steps.
\section{Numerical Results}

To test the effectiveness of the proposed scheme, we have applied it to the simulation of plane wave scattering from a PEC sphere and a space shuttle model. All geometries have been excited by a pulse Gaussian plane wave
\begin{equation}
 \mathbf{e}^{\text{inc}}(\mathbf{r},t)=A_0 \exp \Big(
 \minus\frac{\big(t\minus\frac{\mathbf{\hat{k}}\cdot \mathbf{r}}{c}\big)^2}{2\sigma^2}\Big) \mathbf{\hat{p}}\,,
\end{equation}
where $\sigma=\SI{3620}{\nano\second}$, $\mathbf{\hat{p}}=\mathbf{\hat{x}}$, $\mathbf{\hat{k}}=\minus\mathbf{\hat{z}}$, and $A_0=\SI{1}{\volt\per\meter}$. Three TD-EFIE formulations have been tested: the time-differentiated one, a formulation regularized using quasi-Helmholtz-projectors \cite{dely2019large}, and the new Calder\'on one.

\begin{figure}
\setlength{\abovecaptionskip}{0pt}
\begin{minipage}[c]{\columnwidth}
 \centerline{
 \ref{named}}
 \end{minipage}
 \hfill
 
 \pgfplotsset{small,samples=10}
 \begin{minipage}[c]{.48\columnwidth}
 \centering
 \begin{tikzpicture}
\begin{loglogaxis}[font=\Large,width=\columnwidth,
xtick={10e-9,10e-7,10e-5,10e-3,10e-1},xlabel={$\Delta t $ (\si{\second})},width=\columnwidth,
ylabel={Condition number},
ytick={1,10e5,10e10,10e15,10e19},grid=major
]
\addplot [ultra thick,blue, mark=*,mark size=1.5pt] coordinates {\NsTD};
\addplot [ultra thick,purple, mark=diamond*,mark size=2pt] coordinates {\NsReg};
\addplot [ultra thick,orange, mark=*,mark size=1.5pt] coordinates {\NsCal};
\end{loglogaxis}
%\node[inner sep=0pt] (russell) at (0.60,1.5)
% {\includegraphics[width=30pt]{sphere_IPS_2022.PNG}};

\end{tikzpicture}
\caption{Condition number with respects to the time step ($Ns=270$).}
\label{Cond Dt}
 \end{minipage} \hfill
 \begin{minipage}[c]{.48\columnwidth}
 \centering
 \begin{tikzpicture}
\begin{semilogyaxis}[font=\Large,width=\columnwidth,
xlabel={$h^{-1}$ (\si{\metre})},ylabel={Condition number},ymin=1,
legend columns=-1,
legend entries={Time differentiate,quasi-Helmholtz,This work},
legend to name=named,grid=major
]
\addplot [ultra thick,blue, mark=*,mark size=2pt]
coordinates {\dtTD}; 
\addplot [ultra thick,purple, mark=diamond*,mark size=2pt]
 coordinates {\dtreg};
 \addplot [ultra thick, orange, mark=*,mark size=2pt] 
 coordinates {\dtcal};
\end{semilogyaxis}
%\node[inner sep=0pt] (russell) at (2,1.1)
% {\includegraphics[width=30pt]{sphere_IPS_2022.PNG}};
\end{tikzpicture}
\caption{Condition number with respects to $Ns$ ($\Delta t= \SI{573}{\nano\second}$).}
\label{Cond Ns}
 \end{minipage}
\hfill
 \pgfplotsset{footnotesize,samples=10,legend style={font=\small},width=9cm,height=6cm}
 \begin{minipage}[c]{1\columnwidth}
 \centering
 \begin{tikzpicture}
\begin{semilogyaxis}[xlabel={time (\si{\second})},ylabel={Current intensity $|\bm{j}(t)|$ (\si{\ampere\per\meter})},legend entries={Time differentiated, quasi-Helmholtz, This work},
width=0.98\columnwidth,
legend style={font=\footnotesize,at={(186pt,-5.5pt)}},
ymin=10e-24,ymax=10e-2,ytick={10e-22,10e-15,10e-10,10e-5,10e-2},%grid=major
]
\addplot [ultra thick,blue, mark=*,mark size=3pt] coordinates{(-3.819718634205488e-5, 2.068753634242219e-24)
 (-2.960281941509253e-5, 1.4495780735807236e-24)};
\addplot [ultra thick,purple, mark=diamond*,mark size=3.5pt]coordinates{(-3.819718634205488e-5, 2.068753634242219e-24)
 (-2.960281941509253e-5, 1.4495780735807236e-24)};
\addplot [ultra thick, orange, mark=*,mark size=3pt] coordinates {(-3.819718634205488e-5, 2.068753634242219e-24)
 (-2.960281941509253e-5, 1.4495780735807236e-24)};
\addplot [only marks,blue, mark=*,mark size=3pt] coordinates {\pointTD};
\addplot [only marks,purple, mark=diamond*,mark size=3.5pt] coordinates {\pointReg};
\addplot [only marks, orange, mark=*,mark size=3pt] coordinates {\pointCal};
\addplot [ultra thick,purple] coordinates {\CurrentReg};
\addplot [ultra thick, orange] coordinates {\CurrentCal};
\addplot [ultra thick,blue] coordinates {\CurrentTD};
\end{semilogyaxis}
\node[inner sep=0pt] (russell) at (2.51,1.1)
 {\includegraphics[width=70pt]{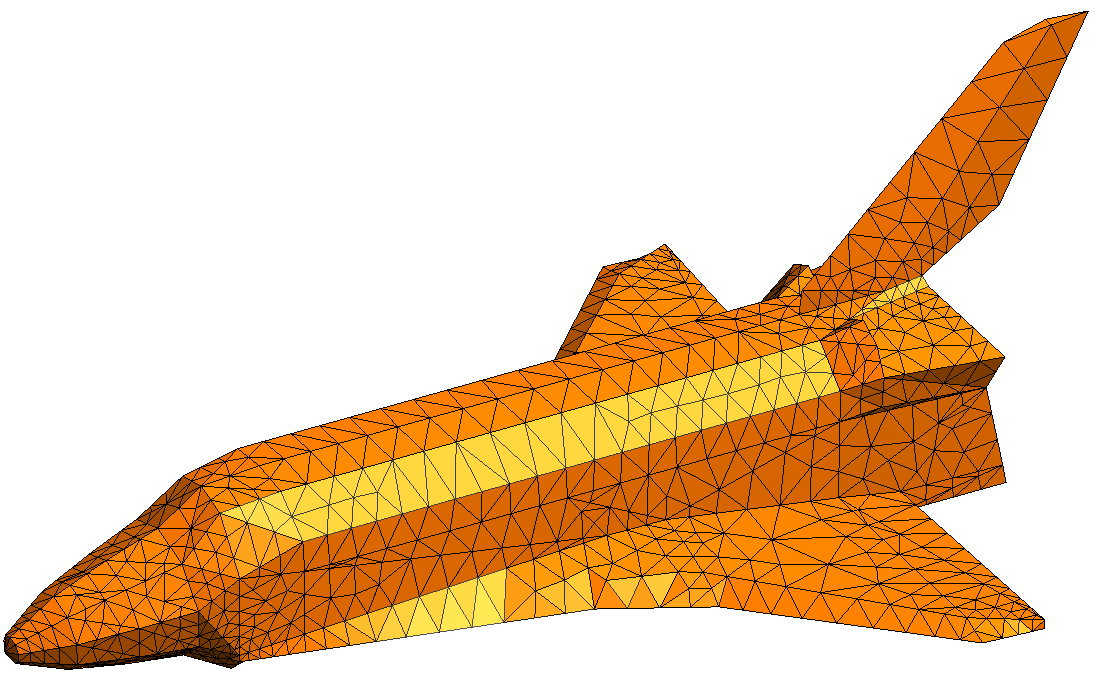}};
\node[ellipse,draw=red, anchor=south,minimum width = 3.5cm, minimum height = 0.7cm,ultra thick] (e) at (5.1,2.3) {};
\node at (5.1,2) {\textcolor{red}{DC-instability}};
\end{tikzpicture}
\caption{Evolution in time of the current intensity at a specific point of the plane with parameters $Ns=4311$ and $\Delta t=\SI{572}{\nano\second}$.}
\label{Current}
\end{minipage}
\vspace{-0.4cm}
\end{figure}

The preconditioning effect of the method we propose has been tested on a spherical scatterer with respect to both the temporal step (Fig.~\ref{Cond Dt}) and mesh refinement (Fig.~\ref{Cond Ns}). These results clearly show that the time-differentiated formulation is the only ill-conditioned one for large time steps. By increasing the refinement of the mesh, however, the conditioning of the quasi-Helmholtz formulation also deteriorates. The Calder\'on TD-EFIE formulation we propose in this work is, therefore, the only one which does not suffer from ill-conditioning due to both large time steps and dense meshes.

To show the favourable properties of our new approach even as pertains DC-instabilities, we have simulated the space shuttle illustrated in Fig.~\ref{Current}. Clearly, while the time differentiated TD-EFIE suffers from DC instabilities, the Calder\'on and regularized formulations we propose is immune from them. Moreover, the the Calder\'on scheme exhibits the best conditioning of \num{41} against \num{1100} for the quasi-Helmholtz and \num{3.5e5} for the time-differentiated ones.

\section*{Acknowledgment}
%\vspace{2.5cm}
The work of this paper has received funding from the EU H2020 research and innovation programme under the Marie Skłodowska-Curie grant agreement n° 955476 (project COMPETE) and from the European Research Council (ERC) under the European Union’s Horizon 2020 research and innovation programme (grant agreement No 724846, project 321).

\bibliographystyle{IEEEtran}
\bibliography{main}

\begin{thebibliography}{1}
\providecommand{\url}[1]{#1}
\csname url@samestyle\endcsname
\providecommand{\newblock}{\relax}
\providecommand{\bibinfo}[2]{#2}
\providecommand{\BIBentrySTDinterwordspacing}{\spaceskip=0pt\relax}
\providecommand{\BIBentryALTinterwordstretchfactor}{4}
\providecommand{\BIBentryALTinterwordspacing}{\spaceskip=\fontdimen2\font plus
\BIBentryALTinterwordstretchfactor\fontdimen3\font minus
  \fontdimen4\font\relax}
\providecommand{\BIBforeignlanguage}[2]{{%
\expandafter\ifx\csname l@#1\endcsname\relax
\typeout{** WARNING: IEEEtran.bst: No hyphenation pattern has been}%
\typeout{** loaded for the language `#1'. Using the pattern for}%
\typeout{** the default language instead.}%
\else
\language=\csname l@#1\endcsname
\fi
#2}}
\providecommand{\BIBdecl}{\relax}
\BIBdecl

\bibitem{lubich1988convolution}
C.~Lubich, ``Convolution quadrature and discretized operational calculus. i,''
  \emph{Numerische Mathematik}, vol.~52, no.~2, pp. 129--145, 1988.

\bibitem{wang2011implicit}
X.~Wang and D.~S. Weile, ``Implicit runge-kutta methods for the discretization
  of time domain integral equations,'' \emph{IEEE transactions on antennas and
  propagation}, vol.~59, no.~12, pp. 4651--4663, 2011.

\bibitem{dely2019large}
A.~D{\'e}ly, F.~P. Andriulli, and K.~Cools, ``Large time step and dc stable
  td-efie discretized with implicit runge--kutta methods,'' \emph{IEEE
  Transactions on Antennas and Propagation}, vol.~68, no.~2, pp. 976--985,
  2019.

\end{thebibliography}
%\begin{thebibliography}{unsrt}
  %\bibitem{dely2019large}
    %D{\'e}ly, A. and Andriulli, F. P and Cools, K., "Large time step and DC stable %TD-EFIE discretized with implicit Runge--Kutta methods" (CRC Press, USA, 1990).
%\end{thebibliography}

\end{document}